\documentclass{amsart}
\usepackage{graphicx}
\usepackage{amsmath}
\usepackage{amsthm}
\usepackage{amsfonts}
\usepackage{dsfont}
\usepackage{amssymb}
\usepackage{microtype}
\usepackage[dvipsnames]{xcolor}
\usepackage{verbatim}
\usepackage{enumitem}
\usepackage{xfrac}
\usepackage{tikz}
\usetikzlibrary{arrows,positioning}
\usepackage{enumitem}
\usepackage{fancyhdr}
\usepackage{lastpage}
\usepackage{amsmath}
\usepackage{amssymb}
\usepackage{amsthm}
\usepackage{multicol}
\usepackage[utf8]{inputenc}
\usepackage{xfrac}
\usepackage{listings}
\usepackage{xcolor}
\usepackage{pgfplots}
\usepackage{etoolbox}
\usepackage{float}
\pgfplotsset{compat=1.17}

\newcommand{\spacewidth}{\the\fontdimen2\font\space}
\newcommand{\nonum}{\hspace{-\spacewidth}}

\makeatletter

\let\@newamstheorem\newtheorem

\newcommand{\@setnumber}[2]{%
  \@namedef{the#1}{#2}%
}

\newcommand{\newcustomtheorem}[2]{%
    \@newamstheorem{#1@inner}{#2}%
    \newenvironment{#1}[1][\nonum]{%
        \@setnumber{#1@inner}{##1}%
        \expandafter\begin{#1@inner}%
    }{%
        \expandafter\end{#1@inner}%
    }%
}

\renewcommand{\newtheorem}{\@ifstar{\newcustomtheorem}{\@newamstheorem}}

\makeatother

\ifdef{\lemma}{ \undef\lemma}{}
\newtheorem{lemma}{Lemma}
\newtheorem{thm}{Theorem}

\newtheorem*{thm*}{Theorem}
\newtheorem*{lemma*}{Lemma}
\newtheorem*{prop*}{Proposition}
\newtheorem*{cor*}{Corollary}

\theoremstyle{definition}
\ifdef{\note}{ \undef\note}{}
\newtheorem*{note}{Note}
\newtheorem*{defn}{Definition}
\newtheorem*{examp}{Example}

\newtheorem*{prob*}{In Class Problem}

\usepackage{dsfont}
\newcommand{\0}{\text{\usefont{U}{bbold}{m}{n}0}}
\newcommand{\1}{\text{\usefont{U}{bbold}{m}{n}1}}

\newcommand{\makemathbb}[1]{%
    \expandafter\newcommand\csname #1\endcsname{\mathbb{#1}}%
}%
\newcommand{\makemathbf}[1]{%
    \expandafter\newcommand\csname b#1\endcsname{\mathbf{#1}}%
}%
\newcommand{\makemathcal}[1]{%
    \expandafter\newcommand\csname c#1\endcsname{\mathcal{#1}}%
}%
\newcommand{\makemathfrak}[1]{%
    \expandafter\newcommand\csname f#1\endcsname{\mathfrak{#1}}%
}%

\let\P\undefined

\makeatletter
\count@=0
\loop
\advance\count@ 1
\edef\letter{\@Alph\count@}%
\expandafter\makemathbb\letter
\expandafter\makemathbf\letter
\expandafter\makemathcal\letter
\expandafter\makemathfrak\letter
\ifnum\count@<26
\repeat
\makeatother

\makemathbf{n}

\DeclareMathOperator*{\argmax}{arg\,max}

\makeatletter
\newtheorem*{rep@theorem}{\rep@title}
\newcommand{\newreptheorem}[2]{%
    \newenvironment{rep#1}[1]{%
        \def\rep@title{#2 \ref{##1}}%
        \begin{rep@theorem}
    }%
    {\end{rep@theorem}}
}
\makeatother
\numberwithin{equation}{section}

\def\N{\mathbb{N}}
\def\P{\mathbb{P}}
\def\R{\mathbb{R}}
\def\Z{\mathbb{Z}}

\def\1{\mathds{1}}
\def\0{\mathbf{0}}

\def\cF{\mathcal{F}}

\def\cN{\mathcal{N}}

\definecolor{darkerred}{RGB}{192,0,0}
\definecolor{darkerblue}{RGB}{0,0,160}
\definecolor{darkgreen}{RGB}{0,160,0}

\newenvironment{PfofTheorem1}[1]
{\par\vskip2\parsep\noindent{\sc Proof of Theorem\ \ref{Thm1}. }}{{\hfill
$\Box$}
\par\vskip2\parsep}

\newenvironment{PfofTheorem2}[2]
{\par\vskip2\parsep\noindent{\sc Proof of Theorem\ \ref{Thm2}. }}{{\hfill
$\Box$}
\par\vskip2\parsep}

\title{Strong recovery in group synchronization}

\author{Bradley Stich}
\address{Bradley Stich\\
Department of Mathematics and Statistics\\
University of North Carolina at Charlotte \\
9201 University City Blvd.\\
Charlotte, NC 28223}
\email{bstich1@uncc.edu}

\begin{document}

\maketitle
\begin{abstract} 
The group synchronization problem is to estimate unknown group elements at the vertices of a graph when given a set of possibly noisy observations of group differences at the edges. We consider the group synchronization problem on finite graphs with size tending to infinity, and we focus on the question of whether the true edge differences can be exactly recovered from the observations (i.e., strong recovery). We prove two main results, one positive and one negative. In the positive direction, we prove that for a sequence of synchronization problems containing the complete digraph along with a relatively well behaved prior distribution and observation kernel, with high probability we can recover the correct edge labeling. Our negative result provides conditions on a sequence of sparse graphs under which it is impossible to recover the correct edge labeling with high probability.
\end{abstract}

\section{Introduction}
Let $G = (V,E)$ be a finite directed graph, and let $\fG$ be a finite group. We will often refer to the group difference between two elements $a,b \in \fG$, which is defined to be the element $a^{-1}b$. 
In the group synchronization problem, vertices of the graph $G$ are labeled by elements of the group $\fG$ according to some probability distribution, and we are given the noisy observations of the group differences along edges of the graph. The task is then to estimate the correct labeling of the vertices up to a global translation by some group element $g \in \fG$, or equivalently, to estimate the correct labeling of the edges. 

For a concrete example of a group synchronization problem, consider the $3$ by $3$ lattice graph in $\Z^2$ together with the group of two elements. More precisely, let $V = \{(i,j): 1 \le i,j \le 3\}$, and let two vertices be connected by an edge whenever they are distance one apart (where the edge is directed from $u$ to $v$ whenever the coordinates of $u$ are less than or equal to the corresponding coordinates of $v$). Let the vertices be labeled uniformly by elements of the group $\Z/2\Z = (\{0,1\}, \oplus)$. For each $u \in V$, let us denote the group element at $u$ by $x(u)$, and for any edge $(u,v) \in E$, the correct group difference at $(u,v)$ is $-x(u) \oplus x(v)$. Since the group we are working with is $\Z/2\Z$, we have $-x(u) \oplus x(v) = x(u) \oplus x(v)$. Suppose we are also given observations of the relative edge differences as follows: for each edge $(u,v) \in E$, the observation of the edge difference is denoted $Y(u,v)$, and the observations are distributed according to a flip probability as follows. With probability $p \in (0,1)$ the observation is incorrect ($Y(u,v) = x(u) \oplus x(v) \oplus 1$), and with probability $1-p$ the observation is correct ($Y(u,v) = x(u) \oplus x(v)$).
In other words, each observation $Y(u,v)$ is assumed to be the true value $x(u) \oplus x(v)$ plus Bernoulli noise. The objective is then to find a way to estimate the correct edge differences. See Figure 1 for a visualization of this example. 

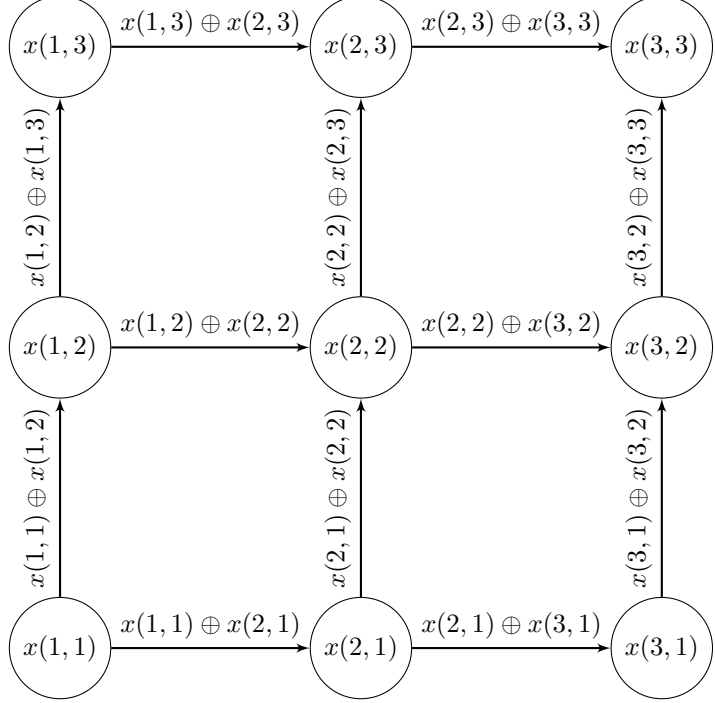
\begin{figure}
    \centering
    
    \caption{Noiseless Synchronization of $\Z/2\Z$ on $3 \times 3 \subset \Z^2$}
    \label{fig:my_label}
\pgfmathsetmacro{\s}{4}
\pgfmathsetmacro{\n}{int(3)}

\begin{center}
    \begin{tikzpicture}
        \tikzstyle{vertex} = [circle,draw,minimum size=1.0em]
        \tikzstyle{edge} = [->,> = latex', thick]
        \tikzstyle{elabel} = [sloped, anchor=center, above]
        \tikzstyle{cont} = [dashed]
        
        \foreach \i in {1, ..., \n} {
            \foreach \j in {1, ..., \n} {
                \pgfmathsetmacro{\x}{\i * \s};
                \pgfmathsetmacro{\y}{\j * \s};
                \node[vertex] (x \i \j) at (\x, \y) {$x(\i,\j)$};
            };
        };
        
        \pgfmathsetmacro{\nmo}{int(\n - 1)}
        \foreach \i in {1, ..., \nmo} {
            \pgfmathsetmacro{\ip}{int(\i + 1)}
            
            \foreach \j in {1, ..., \nmo} {
                \pgfmathsetmacro{\jp}{int(\j + 1)}
                
                \draw[edge] (x \i \j) to node[elabel] {$x(\i,\j) \oplus x(\ip, \j)$} (x \ip \j);
                \draw[edge] (x \i \j) to node[elabel] {$x(\i, \j) \oplus x(\i, \jp)$} (x \i \jp);
            };
            
            \draw[edge] (x \i \n) to node[elabel] {$x(\i, \n) \oplus x(\ip, \n)$} (x \ip \n);
            \draw[edge] (x \n \i) to node[elabel] {$x(\n, \i) \oplus x(\n, \ip)$} (x \n \ip);
        };
    \end{tikzpicture}
\end{center}
\end{figure}

Group synchronization has many applications including network analysis, structure from motion in computer vision, and stochastic physics. We refer the reader to \cite{Cucuringu_2015} for a detailed account of group synchronization in the context of multiplex networks and to \cite{abbe2018group} for connections between group synchronization and stochastic physics. 
Due to its vast applications, group synchronization has received much attention in recent years \cite{el2016graph, ling2020near, dahms2012cluster}. For instance, \cite{lerman2021robust} establishes results on message passing algorithms for solving synchronization problems, and \cite{abbe2020information, polyanskiy2020application} develop and apply information-percolation methods for solving synchronization problems. Additionally, \cite{gao2021geometry} use the classical theory of fibre bundles to characterize the cohomological nature of synchronization. In another recent work, \cite{abbe2018group} asks a question of weak recovery: can the labeling of the vertices be estimated in a manner more accurate than random guessing on the d-dimensional lattice graph, up to a global translation by a group element? 
In contrast, we are interested in the question of strong recovery: given a sequence of synchronization problems, can we recover the vertex labeling up to global translation with probability tending to one? We are particularly interested in determining what effects the geometry of a graph has on our ability to attain strong recovery.

In this work, we prove two main results about strong recovery--one positive and one other negative. For our positive result, we prove that if the graph is the complete digraph on $n$ vertices and the noise is bounded away from a particular constant, then there is a polynomial time algorithm that solves the synchronization problem with probability tending to one as $n$ tends to infinity. In our negative result, we show that if $G_n$ has an independent set $D_n \subset V_n$ whose size tends to infinity, and if the noise does not tend to zero too quickly, then strong recovery is impossible with probability tending to one as $n$ tends to infinity. 

The remainder of the paper is organized as follows. In Section 1.1, we provide definitions for all the necessary objects. In Section 1.2, our main results are presented formally. Section 2 contains the proofs for the positive result, and Section 3 provides the proofs for the negative result.

\subsection{Definitions/Problem Statement}
We will find it convenient to relate a vertex labeling to its corresponding edge labeling. 
Suppose $G = (V,E)$ is a finite connected (directed) graph, with $E \subset V\times V$. Given $x \in \fG^{V}$, define $\psi: \fG^{V} \to \fG^{E}$ by the rule: for each $(u,v) \in E$, we set
\begin{equation*}
    \psi(x)(u,v)=x(u)^{-1}x(v).
\end{equation*}
In what follows, we will assume that $\pi$ is a probability distribution on $\fG$, and we let $\pi^V$ denote the corresponding product distribution on $\fG^V$. We also assume that for each $x \in \fG^V$, we have an observation kernel $Q( \cdot \mid \psi(x))$, which is assumed to be a probability distribution on $\fG^E$. Throughout this paper we use the notation that $X \in \fG^{V}$ is a random variable distributed according to $\pi^{V}$ and $Y \in \fG^{E}$  is a random variable distributed according to $Q(\cdot \mid \psi(X))$.

\begin{defn}[1.1]
An \textit{estimator} $T: \fG^{E} \to \fG^{E}$ is a function which maps the observations onto a configuration of group elements on $E$. 

\end{defn}
An example of an estimator is the trivial estimator which maps the set of observations onto itself. 
For convenience, we now define an object that contains all the relevant information for a synchronization problem.
\begin{defn}[1.2]
A \textit{synchronization problem} $\fQ$ is defined as a tuple
$$\fQ = (G, \fG, Q, \pi),$$
where $G = (V,E)$ is a graph, $\fG$ is a finite group, $Q$ is an observation kernel, and $\pi$ is a prior distribution on $\fG$.
\end{defn}
In this work we are concerned with sequences of synchronization problems, $\{\fQ_n\}_{n=1}^{\infty} = \{(G_n,\fG,Q_n,\pi_n)\}_{n=1}^{\infty}$, in which the group $\fG$ is fixed. 
\begin{defn}[1.3]
Let $\{\fQ_n\} = \{(G_n,\fG,Q_n,\pi_n)\}$ be a sequence of synchronization problems. For each $n \in \N$, draw $X_n$ according to $\pi_n^{V_n}$, draw $Y_n$ according to $Q_n(\cdot \mid \psi_n(X_n))$, and let $\theta_n = \psi_n(X_n)$. We say that the sequence of estimators $\{T_n\}_{n=1}^\infty$ has the \textit{strong recovery property} if 
\begin{equation*}
    \lim_{n \to \infty} \P(T_n(Y_n) = \theta_n) \to 1.
\end{equation*}
\end{defn}
This definition states that a sequence of estimators $\{T_n\}$ has the strong recovery property if the probability that the estimator yields the correct configuration grows arbitrarily close to 1 as $n$ grows large. 
\begin{defn}[1.4]
Let $\{\fQ_n\}_{n=1}^{\infty} = \{(G_n, \fG, Q_n, \pi_n)\}_{n=1}^{\infty}$ be a sequence of synchronization problems. We say that the sequence of synchronization problems $\{\fQ_n\}_{n=1}^{\infty}$ has the \textit{strong recovery property} if there exists a sequence of estimators $\{T_n\}_{n=1}^{\infty}$ on $\{\fQ_n\}_{n=1}^{\infty}$ that has the strong recovery property. 
\end{defn}
Essentially, given a sequence of synchronization problems $\{\fQ_n\}$, 
the strong recovery problem is to define a sequence of estimators $\{T_n\}$ which approximate the correct edge configuration with high probability in the limit. 
\begin{defn}[1.5]
We say that $Q$ is a \textit{uniform observation kernel with flip probability $p \in (0,1)$} if the following holds. 
For each edge $(u,v) \in E$, first define $Q'$ as follows:
\begin{equation*}
    Q'(Y(u,v)\mid\psi(x)(u,v)) = \begin{cases}
    1-p, & \text{if } Y(u,v) = x(u)^{-1}x(v) \\
    \frac{p}{|\fG|-1}, & \text{if } Y(u,v) \neq x(u)^{-1}x(v).
                                \end{cases}
\end{equation*}
Then $Q( \cdot \mid \psi(x))$ is a product of these values over all of the edges:
\begin{equation*}
    Q(Y\mid\psi(x)) = \prod_{(u,v) \in E} Q'(Y(u,v)\mid\psi(x)(u,v)).
\end{equation*}
Finally, we say that $Q$ is a uniform observation kernel with flip probability $p$ if $Q( \cdot \mid \psi(x))$ satisfies this property for all $x \in \fG^V$.
\end{defn}
Note that if $Q( \cdot \mid \psi(x))$ is a uniform observation kernel, then the edge observations are all independent of each other (conditional on $\psi(x)$), and for each $(u,v) \in E$, the observation $Y(u,v)$ may be written as 
\begin{equation*}
    Y(u,v) = x(u)^{-1}x(v)\fN_{u,v},
\end{equation*}
where the noise variable $\fN_{u,v} \in \fG$ is equal to the identity element $e$ with probability $1-p$ and with probability $p$ it is uniform on $\fG \setminus \{e\}$. 
A uniform observation kernel will select the correct group difference with probability $1-p$, and it will uniformly select an incorrect group difference otherwise. For instance, if $\fG$ is the group of two elements, then $p$ is just the flip probability. 
\subsection{Main Results}
Our first result provides conditions for which a sequence of synchronization problems has the strong recovery property. For notation, we let $K_n$ be the complete digraph on $n$ vertices, which has vertex set $V_n = \{1,\dots,n\}$ and edge set $E_n = V_n \times V_n$.
\begin{thm} \label{Thm1}
Let $\{\fQ_n\}_{n=1}^{\infty} = \{(K_n, \fG_n, Q_n, \pi_n)\}_{n=1}^{\infty}$, where $K_n$ is the complete digraph with $n$ vertices, $\fG$ is an arbitrary finite group with $|\fG| \ge 2$, for each $n \in \N$, $Q_n$ is a uniform observation kernel with flip probability $p_n$, and $\pi_n$ is the uniform distribution on $\fG$. Further suppose that there exists $\beta > 0$ such that $\{p_n\}_{n=1}^{\infty} \subset [0,1] \setminus (p_c - \beta, p_c + \beta)$, where $p_c := 1 - 1/|\fG|$. Then $\{\fQ_n\}$ has the strong recovery property. 
\end{thm}

This result indicates that the complete graph provides enough redundant information that strong recovery is possible as long as the noise is bounded away from the critical constant $p_c$.

Our second result provides conditions for which the a sequence of synchronization problems cannot have the strong recovery property. 
We require a few definitions before stating the result. For a directed graph $G = (V,E)$ and a vertex $v \in V$, the degree of $v$, denoted $\text{deg}(v)$, is the number of edges incident to $v$. We say that a sequence of graphs $\{G_n\}_{n=1}^{\infty}$ has bounded degree if there exists $d \in \N$ such that for all $n \in \N$ and all $v \in V_n$, $\text{deg}(v) \le d$, and say that $V_n$ is bounded by $d$. A set $D \subset V$ is said to be independent if for each pair of vertices $u,v \in D$, there is no edge $(u,v) \in E$. Also, for non-negative sequences  $\{a_n\}_{n=1}^{\infty}$ and $\{b_n\}_{n=1}^{\infty}$,  we say that $b_n = \omega(a_n)$ if $\forall c\in \R$ with $ c>0$, $\exists N\in \N$ such that $\forall n \ge N, b_n > ca_n \ge 0$. 
\begin{thm} \label{Thm2}
Let $\{\fQ_n\}_{n=1}^{\infty} = \{(G_n, \fG, Q_n, \pi_n)\}_{n=1}^{\infty}$ be a sequence of synchronization problems where  $G_n = (V_n,E_n)$ is a simple directed graph which has bounded degree with $V_n$ bounded by $d \in \N$, $\fG$ is a finite group with at least 2 elements, $\pi_n$ is the uniform distribution on $\fG$ and $Q_n$ is the uniform observation kernel with flip probability $p_n$. If
\begin{itemize}
    \item there exists an independent set $D_n \subset V_n$ such that $|D_n| \to \infty$ as $n \to \infty$,
    \item $\{p_n\}_{n=1}^\infty \subset (0, \sfrac{1}{2})$, and
    \item $p_n = \omega\bigl(|D_n|^{-1/d}\bigr)$,
\end{itemize}
then $\{\fQ_n\}$ does not have the strong recovery property. 
\end{thm}

This result tells us that for bounded degree sequences, even if the noise converges to zero, provided that it does so sufficiently slowly, $\{\fQ_n\}$ cannot have the strong recovery property due to excessive error propagation. 
\begin{examp}
Let $G_n \subset \Z^3$ be the discrete hypercube with side length $n$, which has vertex set $V_n = \{(i,j): 1 \le i,j \le n\}$ and edge set $E_n$ in which any two vertices are connected by an arbitrarily oriented edge whenever they are distance one apart. Consider $\{\fQ_n\} = \{(G_n, S_5, Q_n, \pi_n)\}$, where $S_5$ is the group of permutations on a set of five elements, $Q_n$ is a uniform observation kernel and $\pi_n$ is the uniform distribution on $S_5$. First note that $d$ is bounded by $6$. Also, for each $n \in \N$, it is possible to construct an independent subset $D_n \subset V_n$ such that $|D_n| \approx kn^3$ for some constant $k \in \R$. Therefore, if $\{p_n\}_{n=1}^\infty \subset (0, \sfrac{1}{2})$ and $p_n = \omega\bigl(n^{3(-1/6)}\bigr) =  \omega\bigl(n^{-1/2}\bigr)$, then $\{\fQ_n\}$ does not have the strong recovery property by Theorem 2. 
\end{examp}
\section{Synchronization on Complete Graphs} 
The basic idea of the proof of Theorem 1 is to define a sequence of ``triangle estimators" on $\{\fQ_n\}$, and show that it has the Strong Recovery Property. Note that for any $(u,v) \in E_n$, and for any vertex $w \in V_n\setminus \{u,v\}$, in the case of noiseless observations we have
\begin{equation*}
    Y_n(u,w) Y_n(w,v) = x_n(u)^{-1} x_n(w) x_n(w)^{-1} x_n(v) = x_n(u)^{-1}x_n(v) = Y_n(u,v).
\end{equation*}
Thus, we can consider all products of observations of the form $Y_n(u,w) Y_n(w,v)$ and base our estimate for the group difference associated to edge $(u,v)$ off of them, rather than the observation $Y_n(u,v)$. In doing so, we achieve a collection of ``votes", and we will show that with high probability these votes provide good estimates of the true group differences.
For each edge $(u,v) \in E_n$, the triangle estimator $T_n$ counts for each $g\in \fG$ the number of vertices $x\in V_n$ such that $Y_n(u,x)Y_n(x,v) = g$ and assigns the edge $(u,v)$ the group difference with the most votes. Then we show that as we take the limit as $n$ goes to infinity, the probability that the triangle estimator is incorrect on any edge goes to zero. 

For the argument showing the efficacy of the sequence of triangle estimators, we need to find the probability that an arbitrary pair of observations between the vertices $1,2 \in V_n$ is either correct, or yields a particular incorrect element.  
For the following lemma, assume the hypotheses of Theorem 1. 
\begin{lemma} \label{Lemma1} Consider the edge $(1,2) \in E_n$ and let $v\in V_n \setminus \{1,2\}$. Let $g^* := X_n(1)^{-1}X_n(2)$, and let $g \in \fG \setminus \{g^*\}$ Then: 
\begin{align*}
    & \P(Y_n(1,v)Y_n(v,2) = g^*) = 1-2p_n + p_n^2\bigg(\frac{|\fG|}{|\fG|-1} \bigg), \text{ and }\\
    & \P(Y_n(1,v)Y_n(v,2) = g) =  \frac{2(p_n-p_n^2)}{|\fG|-1} + \frac{p_n^2(|\fG| - 2)}{(|\fG|-1)^2}.
\end{align*}

\end{lemma}

\begin{proof}
Let $X_n(1)^{-1}X_n(v) = g_{1v}$, and $X_n(v)^{-1}X_n(2) = g_{v2}$, so that $g^* = g_{1v}g_{v2}$. We define the following three events:
\begin{align*}
    & S := \{Y_n(1,v)Y_n(v,2) = g^*\}\\
    & \fA := \{Y_n(1,v) = g_{1v}, Y_n(v,2) = g_{v2}\}\\
    & \fB := \{Y_n(1,v)Y_n(v,2) = g^* \} \cap \{Y_n(1,v) \ne g_{1v}, Y_n(v,2) \ne g_{v2}\}.
\end{align*}
Then $S$ can be expressed as the union of disjoint events:
\begin{equation*}
    S = \fA \cup \fB \implies \P(S) = \P(\fA) + \P(\fB). 
\end{equation*}
By the independence of observations, the probability of $\fA$ is precisely $(1-p_n)^2$. 

The probability of $\fB$ is given by the probability of choosing a pair $(g',g'') \in (\fG \setminus \{g_{1v}\}) \times (\fG \setminus \{g_{v2}\})$ such that $g'g'' = g^*$. 

If we take any particular $g' \in \fG \setminus \{g_{1v}\}$, and $g'' \in \fG \setminus \{g_{v2}\}$, then by the independence of observations, 

\begin{equation*}
    \P(Y_n(1,v) = g', Y_n(v,2) = g'') = \bigg( \frac{p_n}{|\fG| - 1} \bigg)^2.
\end{equation*}

If we fix $g' \in \fG \setminus \{g_{1v}\}$, then there is a unique $g'' \in \fG \setminus \{g_{v2}\}$, such that $g'g'' = g^*$ due to cancellation. There are $|\fG| - 1$ choices for $g'$, therefore, there are $|\fG| - 1$ such pairs, and we see that 

\begin{equation*}
    \P(\fB) = \bigg( \frac{p_n}{|\fG| - 1} \bigg)^2 (|\fG| - 1) = \frac{p_n^2}{|\fG| - 1}.
\end{equation*}
Therefore, 
\begin{equation*}
    \P(S) = (1-p_n)^2 + \frac{p_n^2}{|\fG| - 1} = 1-2p_n + p_n^2\bigg(\frac{|\fG|}{|\fG|-1} \bigg).
\end{equation*}

When $g \in \fG \setminus \{g^*\}$, we can find $\P(Y_n(1,v)Y_n(v,2) = g)$ by similar direct calculation, or by noting that the compliment of the event $S = \{Y_n(1,v)Y_n(v,2) = g^*\}$ is the event that $\{Y_n(1,v)Y_n(v,2) = g_0\}$ for some $g_0 \in \fG \setminus \{g^*\}$, and this distribution is uniform on $\fG \setminus \{g^*\}$, so using algebra, it can be shown that
\begin{equation*}
    \P(Y_n(1,v)Y_n(v,2) = g) = \P(S^c) / (|\fG| - 1) = \frac{1 - \P(S)}{|\fG| - 1} = \frac{2(p_n-p_n^2)}{|\fG|-1} + \frac{p_n^2(|\fG| -2)}{(|\fG|-1)^2}.  
\end{equation*}
\end{proof}

\begin{lemma} \label{Lemma2}
       Let $\delta>0$, let $p_c := 1 - 1/|\fG|$ and let $f,h:\R \to \R$ be defined 
       \begin{align*}
           & f(x) = 1-2x + x^2\bigg(\frac{|\fG|}{|\fG|-1} \bigg), \\
           & h(x) = \frac{2(x-x^2)}{|\fG|-1} + \frac{x^2(|\fG| - 2)}{(|\fG|-1)^2}.
       \end{align*}
           Then whenever $x \in \R$ with $x \notin (p_c - \delta, p_c + \delta)$, the following holds: 
           \begin{align*}
               & f(x) \ge \min(f(p_c - \delta), f(p_c + \delta)) > f(p_c) = 1/|\fG|,\\
               & h(x) \le \max(h(p_c - \delta), h(p_c + \delta)) < h(p_c) = 1/|\fG|. 
           \end{align*}
       
\end{lemma}

\begin{proof}
Note that $f$ and $h$ are polynomials, and therefore, are continuous and differentiable on $\R$. We find 
\begin{align*}
    & f'(x) = 2x \frac{|\fG|}{|\fG| - 1} - 2, \\
    & h'(x) = \frac{2-4x}{|\fG| -1} + \frac{2x(|\fG|-2)}{(|\fG|-1)^2}.
\end{align*}

Since $f,h$ are quadratic functions and $f'(p_c) = 0 = h'(p_c)$, we see that both $f$ and $h$ have their extremal value at $p_c$. Now, $f'(x) < 0$ for $x \in (0, p_c)$ and $f'(x) > 0$ for $x \in (p_c, 1)$, so $f$ is strictly decreasing on $(0, p_c)$ and strictly increasing on $(p_c,1)$. A similar check shows that $h$ is strictly increasing on $(0,p_c)$ and strictly decreasing on $(p_c,1)$. Hence, $f(p_c \pm \delta) > f(p_c)$, and $f(x) \ge \min(f(p_c - \delta), f(p_c + \delta))$ when $x \notin (p_c - \delta, p_c + \delta)$. Similarly, $h(x) < h(p_c)$ and $h(x) \le \max(h(p_c - \delta), h(p_c + \delta)) < h(p_c) = 1/|\fG|$ for $x \notin (p_c - \delta, p_c + \delta)$.

\end{proof}

\begin{PfofTheorem1}

Assume the hypotheses of Theorem 1. First we must define an estimator on $\{\fQ_n\}$. A 2-observation on an edge $(u,v) \in E_n$ is a pair of observations of the form $(Y_n(u,w), Y_n(w,v))$. Note that if the 2-observation $(Y_n(u,w), Y_n(w,v))$ on the edge $(u,v)$ is correct, then $Y_n(u,w)Y_n(w,v) = X_n(u)^{-1}X_n(v)$. Let $T_n: \fG^{E_n}\to\fG^{E_n}$ be defined as the estimator which chooses the configuration such that each edge is chosen as the group difference given by the most 2-observations. Let $(u,v) \in E_n$. Omitting the indices for brevity, for each integer $n \ge 3$ and for each $g \in \fG$ we define the function $\gamma_{u,v}^{(g)}:V_n\setminus\{u,v\} \to \{0,1\} \subset \Z$ to be:
\begin{equation*}
    \gamma_{u,v}^{(g)}(w) = \begin{cases} 
      1 & \text{if } Y_n(u,w)Y_n(w,v) = g \\
      0 & \text{otherwise}. 
   \end{cases}
\end{equation*}
Again omitting the indices for brevity, for each integer $n \ge 3$ let $\phi_{u,v}:\fG \to \{0,\dots,n-2\} \subset \Z$ be defined
\begin{equation*}
    \phi_{u,v}(g) = \sum_{w \in V_n \setminus \{u,v\}} \gamma_{u,v}^{(g)}(w).
\end{equation*}
For each $(u,v) \in E_n$,
\begin{equation*}
    T_n(Y_n)(u,v) = \argmax_{g \in \fG} \phi_{u,v}(g),
\end{equation*}
with ties broken arbitrarily.
Let $F_n(u,v)$ be the event that the estimator is incorrect on the edge $(u,v) \in E_n$. That is, 
\begin{equation*}
    F_n(u,v) = \{T_n(Y_n)(u,v) \ne X_n(u)^{-1}X_n(v)\}.
\end{equation*}
By the symmetry of $K_n$, we note that $\P(F_n(u,v)) = \P(F_n(1,2))$ for any edge $(u,v) \in E_n$.
Let $A_n$ be the event that the estimator is incorrect on any edge of $E_n$:
\begin{equation*}
    A_n = \bigcup_{(u,v) \in E_n}F_n(u,v).
\end{equation*}
Since there are $\binom{n}{2}$ edges in the undirected complete graph, for the complete digraph $|E_n| = 2 \binom{n}{2}$. By the union bound, we have:
\begin{equation*}
    \P(A_n) \le \sum_{(u,v) \in E_n}\P(F_n(u,v)) = 2 \binom{n}{2}\P(F_n(1,2)) = n(n-1)\P(F_n(1,2)).
\end{equation*}

For any $g \in \fG$, we define the event 

\begin{equation*}
    W_g = \bigg\{ \phi_{1,2}(g) \le \frac{n-2}{|\fG|} \bigg\}.
\end{equation*}
Define $g^* = X_n(1)^{-1}X_n(2)$. 
Note that if the triangle estimator gives more than $(n-2)/|\fG|$ votes to $g^*$ and less than or equal to $(n-2)/|\fG|$ votes to any other element of $\fG$, then the triangle estimator is correct on the edge $(1,2)$:
\begin{equation*}
            W_{g^*}^c \cap \left( \bigcap_{g \neq g^*} W_g \right) \subset F_n(1,2)^c.
\end{equation*}
Taking compliments, we see that
\begin{equation*}
        F_n(1,2) \subset W_{g^*} \cup \left( \bigcup_{g \neq g^*} W_g^c \right). 
\end{equation*}
By the union bound, we have
\begin{equation*}
    \P(F_n(1,2)) \le \P(W_{g^*}) + \sum_{g \ne g^*} \P(W_g^c). 
\end{equation*}

Recall that by Lemma 1, the probability that the triangle estimator is correct on the edge $(1,2)$ is given by the quadratic polynomial in $p_n$:

\begin{equation*}
    f(p_n) := 1-2p_n + p_n^2\bigg(\frac{|\fG|}{|\fG|-1} \bigg).
\end{equation*}
Similarly, the probability that the triangle estimator yields any particular $g \in \fG \setminus\{g^*\}$ is given by the quadratic polynomial in $p_n$:

\begin{equation*}
    h(p_n) := 2 \frac{(p_n-p_n^2)}{|\fG|-1} + \frac{p_n^2}{(|\fG|-1)^2}.
\end{equation*}

Note that $\phi_{1,2}(g^*)$ is the sum of the $n-2$ i.i.d. Bernoulli random variables $\gamma_{1,2}^{g^*}(v)$ (one for each vertex $v \in V^n \setminus \{1,2\}$), and that $\gamma_{1,2}^{g^*}(v) = 1$ whenever the $2$-observation $1$-$v$-$2$ is correct, which happens with probability $f(p_n)$. Let $\epsilon_n = f(p_n) - 1/|\fG|$. Then since $p_n \notin (p_c - \beta, p_c + \beta)$, by Lemma 2, there exists $\lambda_1$ such that $\epsilon_n > \lambda_1 > 0$, so that Hoeffding's inequality gives us: 
\begin{equation*}
    \P(W_{g^*}) = \P\bigg(\phi_{1,2}(g^*) \le \frac{n-2}{|\fG|}\bigg) \le 2e^{-2\epsilon^2(n-2)} \le 2e^{-2\lambda_1^2(n-2)}.
\end{equation*}

Let $g \in \fG\setminus\{g^*\}$. Then the following holds:
\begin{equation*}
    W_g^c = \bigg\{\phi_{1,2}(g) > \frac{n-2}{|\fG|}\bigg\} \subset \bigg\{\phi_{1,2}(g) \ge \frac{n-2}{|\fG|}\bigg\} := W_g^{c'}. 
\end{equation*}

Note that $\phi_{1,2}(g)$ is the sum of $n-2$ i.d.d. Bernoulli random variables $\gamma_{1,2}^g(v)$ and that $\gamma_{1,2}^g(v) = 1$, which happens with probability $h(p_n)$. Let $\hat{\epsilon}_n = 1/|\fG| - h(p_n)$. Then since $p_n \notin (p_c - \beta, p_c + \beta)$, by Lemma 2, there exists $\lambda_2$ such that $\hat{\epsilon}_n > \lambda_2 > 0$, and by Hoeffding's inequality we have: 

\begin{equation*}
    \P(W_{g}^c) \le \P(W_g^{c'}) = \P\bigg(\phi_{1,2}(g) \ge \frac{n-2}{|\fG|}\bigg) \le 2e^{-2\hat{\epsilon}_n^2(n-2)} \le 2e^{-2\lambda_2^2(n-2)}.
\end{equation*}

Set $\lambda := \min(\lambda_1^2, \lambda_2^2)$. Then the following holds:

\begin{align*}
    \P(A_n) & \le n(n-1)\P(F_n(1,2))\\
            & \le n(n-1)\bigg(\P(W_{g^*}) + \sum_{g \ne g^*}\P(W_g^c)\bigg) \\
            & \le n(n-1)\bigg(\P(W_{g^*}) + \sum_{g \ne g^*}\P(W_g^{c'})\bigg)\\
            & \le 2n(n-1)(e^{-2\lambda_1^2(n-2)} + (|\fG|-1)e^{-2\lambda_2^2(n-2)}) \\
            & \le 2n(n-1)(e^{-2\lambda(n-2)} + (|\fG|-1)e^{-2\lambda(n-2)}) \\
            & = \frac{2n(n-1)|\fG|}{e^{2\lambda(n-2)}}. 
\end{align*}

Taking L'H$\hat{o}$pital's rule twice, we find that 

\begin{equation*}
    \lim_{n \to \infty} \P(A_n) \le \lim_{n \to \infty} \frac{2n(n-1)|\fG|}{e^{2\lambda(n-2)}} = \lim_{n \to \infty} \frac{|\fG|}{\lambda^2 e^{2\lambda(n-2)}} = 0.
\end{equation*}
    
Therefore as $n$ tends to infinity $\P(A_n)$ tends to zero, and we conclude that the probability that the estimator $T_n$ is correct for every edge $(u,v) \in E_n$ goes to one as n goes to infinity:  
\begin{equation*}
    \lim_{n \to \infty} \P(T_n(Y_n) = \psi(X_n)) \to 1.
\end{equation*}
As a result, the sequence $\{T_n\}$ has the strong recovery property and we conclude that $\{\fQ_n\}$ has the strong recovery property.

\end{PfofTheorem1}

\section{Impossibility of Synchronization on Sparse Graphs}

For each of the following lemmas, assume the hypotheses of Theorem 2. Our argument will show that the \textit{maximum a posteriori} (MAP) estimator fails with near certainty as $n$ grows large. We first define the posterior distribution on $\fG^{V_n}$ as it will be necessary to reference as we derive the MAP estimator. Recall that $Q_n$ is a uniform observation kernel. Given $y_n \in \fG^{E_n}$, the posterior distribution on $\fG^{V_n}$ is given as follows: $\forall x_n \in \fG^{V_n}$, 
\begin{equation*}
    \pi_n'(x_n|y_n) = \frac{Q_n(y_n\mid\psi(x_n))\pi_n(x_n)}{\sum_{x_n'\in \fG^{V_n}}Q_n(y_n|\psi(x_n'))\pi_n(x_n')}.
\end{equation*}
The maximum a posterior estimator on $\fG^{V_n}$ is defined by:
\begin{equation*}
    \Phi_n(y_n) = \argmax_{x_n' \in \fG^{V_n}} \pi_n'(x_n'\mid y_n),
\end{equation*}
with ties broken arbitrarily.
\begin{defn}[2.1]
Let $x_n, x'_n \in \fG^{V_n}$, let $g \in \fG$. We write $x'_n = g x_n$ if for each $u \in V_n$, $x'_n(u) = gx_n(u)$.
\end{defn}
\begin{lemma}
Let $x_n, x'_n \in \fG^{V_n}$. Then
\begin{equation*}
    \psi(x_n) = \psi(x_n') \iff \exists \: g \in \fG : x_n' = gx_n.
\end{equation*}
\end{lemma}
\begin{proof}
Suppose $x_n'=gx_n$. Then for each $(u,v) \in E_n$, \begin{equation*}
    \psi(x_n')(u,v) = (x_n(u)')^{-1}x_n(v)' = (gx_n(u))^{-1}gx_n(v) =
    \end{equation*}
    \begin{equation*}
    x_n(u)^{-1}(g^{-1}g)x_n(v) = x_n(u)^{-1}x_n(v) = \psi(x_n)(u,v).
\end{equation*}
We leave the other implication as an exercise for the reader. 
\end{proof}
We define an equivalence class for vertex configurations: 
\begin{defn}[2.2]
Let $[x_n] = \{x_n':\exists \: g \in \fG: x_n' = gx_n\}$. 
\end{defn}
\begin{lemma}
If $w_n \in [x_n]$, then $\forall y_n \in \fG^{E_n}$, 
\begin{equation*}
    \pi_n'(w_n\mid y_n) = \pi_n'(x_n\mid y_n).
\end{equation*}
\end{lemma}
\begin{proof}
Since $w_n \in [x_n]$, there exists $g \in \fG$ such that $w_n = gx_n$. Then by Proposition 1, this implies that $\psi(w_n) = \psi(x_n)$. By this result and the fact that $\pi_n$ is a uniform distribution, it follows that:
\begin{align*}
    \pi_n'(x_n\mid y_n) & = \frac{Q_n(y_n\mid \psi(x_n))\pi_n(x_n)}{\sum_{x_n'\in \fG^{V_n}}Q_n(y_n\mid \psi(x_n'))\pi_n(x_n')}\\ & = \frac{Q_n(y_n\mid \psi(w_n))\pi_n(w_n)}{\sum_{x_n'\in \fG^{V_n}}Q_n(y_n\mid \psi(x_n'))\pi_n(x_n')}\\ & = \pi_n'(w_n\mid y_n).
\end{align*}
\end{proof}
We have shown that given $x_n \in \fG^{V_n}$, the MAP estimator will map all $w_n \in [x_n]$ to the same value.     
    \begin{lemma}
        Let $\{p_n\}_{n=1}^\infty \subset (0, \sfrac{1}{2})$, let $\{a_n\}_{n=1}^\infty$ be in $(0,\infty) \subset \R$ such that $a_n \to \infty$, and let $d > 0$. If $p_n = \omega\bigl(a_n^{-1/d}\bigr)$, then for any constant $K > \sfrac{1}{2}$,
        \begin{equation*}
            \lim_{n \to \infty} \left(1 - \left(\frac{p_n}{K}\right)^d\right)^{a_n} = 0.
        \end{equation*}
        \begin{proof}
            First, using that $e^x \ge 1 + x$ (which follows from the fact that $e^x$ is convex, and $1+x$ is the line tangent to $e^x$ at $x = 0$), it follows that
            \begin{equation*}
                e^{rx} = (e^x)^r \ge (1+x)^r.
            \end{equation*}
            Then
            \begin{equation*}
                \left(1 - \left(\frac{p_n}{K}\right)^d\right)^{a_n} \le e^{-a_n p_n^dK^{-d}} = \left(e^{K^{-d}}\right)^{-a_n p_n^d}.
            \end{equation*}
            By assumption, $p_n = \omega\bigl(a_n^{-1/d}\bigr)$, which implies that
            \begin{equation*}
                \lim_{n \to \infty} p_n a_n^{1/d} = \lim_{n \to \infty} \frac{p_n}{a_n^{-1/d}} = \infty,
            \end{equation*}
            and therefore since $d > 0$ is a constant, and $x^d$ is monotonically increasing,
            \begin{equation*}
                \lim_{n \to \infty} p_n^d a_n = \lim_{n \to \infty} \left(p_n a_n^{1/d}\right)^d = \infty.
            \end{equation*}
            Since $e$, $K$, and $d$ are constants, $e^{K^{-d}}$ is also a constant, and since $K > 0$, $K^{-d} > 0$ as well. Since $e^x$ is strictly monotonically increasing, this implies $e^{K^{-d}} > e^0 = 1$. This, along with the fact that $a_n p_n^d \to \infty$, gives that
            \begin{equation*}
                \lim_{n \to \infty} \left(e^{K^{-d}}\right)^{-a_n p_n^d} = 0.
            \end{equation*}
            Therefore,
            \begin{equation*}
                \lim_{n \to \infty} \left(1 - \left(\frac{p_n}{K}\right)^d\right)^{a_n} \le \lim_{n \to \infty} \left(e^{K^{-d}}\right)^{-a_n p_n^d} = 0.
            \end{equation*}
            Since $K > \sfrac{1}{2}$, and $p_n < \sfrac{1}{2}$ for all $n$, $\sfrac{p_n}{K} < 1$ for all $n$. This additionally implies that $(\sfrac{p_n}{K})^d < 1$, which further implies
            \begin{equation*}
                1 - \left(\frac{p_n}{K}\right)^{d} > 0,
            \end{equation*}
            and ultimately
            \begin{equation*}
                \left(1 - \left(\frac{p_n}{K}\right)^d\right)^{a_n} > 0.
            \end{equation*}
            As such,
            \begin{equation*}
                0 \le \lim_{n \to \infty} \left(1 - \left(\frac{p_n}{K}\right)^d\right)^{a_n} \le 0,
            \end{equation*}
            and therefore,
            \begin{equation*}
                \lim_{n \to \infty} \left(1 - \left(\frac{p_n}{K}\right)^d\right)^{a_n} = 0.
            \end{equation*}
        \end{proof}
    \end{lemma}

\begin{lemma}
       Fix $g \in \fG\setminus\{e\}$. As $n$ goes to infinity, with probability tending to one, there exists a vertex $v \in V_n$ such that for all edges of the form $(u,v)$ (or $(v,u)$) in $E_n$, we have $Y_n(u,v) = X_n(u)^{-1}X_n(v)g$ (or $Y_n(v,u) = (X_n(v)g)^{-1}X_n(u)$, respectively).
\end{lemma}   
 In other words, with probability one there exists a vertex such that each observation misidentifies it to be a particular group element. We say that the observations of such a vertex are offset by $g$.  
\begin{proof}
For each $u \in V_n$, let
\begin{equation*}
    A_n^{(g)}(w,u) := \{Y_n(w,u) = X_n(w)^{-1}X_n(u)g\} \cup  \{Y_n(u,w) = X_n(u)^{-1}X_n(w)m_n(w)\},
\end{equation*}
where for each difference $X_n(u)^{-1}X_n(w)$, we define 
\begin{equation*}
m_n(w) := X_n(w)^{-1}X_n(u)g^{-1}X_n(u)^{-1}X_n(w)
\end{equation*}
to account for the different edge orientations. Note that $X_n(u)^{-1}X_n(w)m_n(w) = g^{-1}X_n(u)^{-1}X_n(w) = (X_n(u)g)^{-1}X_n(w)$, and that if we assume that the graph $G_n$ is undirected, with the assumption that $Y_n(u,v) = Y_n(v,u)^{-1}$, then the proof of this lemma simplifies as the $m_n$ term is no longer necessary. Also, observe that 
\begin{equation*}
    \P(A_n^{(g)}(w,u)) = \frac{p_n}{|\fG|-1}.
\end{equation*}
We will use the following notation for the neighborhood of a vertex $u$: $N_n(u) = \{v \in V_n : (v,u) \in E_n \text{ or } (u,v) \in E_n \}$.
Let 
\begin{equation*}
    B_n^{(g)}(u) := \bigcap_{w\in N_n(u)}A_n^{(g)}(w,u).
\end{equation*}
Then since the noise is independent on edges and all degrees are bounded by $d$, we have
\begin{equation*}
\P(B_n^{(g)}(u)) = \left(\frac{p_n}{|\fG|-1}\right)^{deg(u)} \ge \left(\frac{p_n}{|\fG|-1}\right)^{d} >0.
\end{equation*}
The probability that all of the $deg(u)$ observations on the edges $(u,w)$ are not either $X_n(w)^{-1}X_n(u)g$ or $(X_n(u)g)^{-1}X_n(w)$ is 
\begin{equation*}
\P(B_n^{(g)}(u)^c) = 1-\P(B_n^{(g)}(u)) = 1-\left(\frac{p_n}{|\fG|-1}\right)^{deg(u)} < 1-\left(\frac{p_n}{|\fG|-1}\right)^d < 1.
\end{equation*}
Let $\cF_n^{(g)}$ be the event that there exists a $u \in D_n$ such that $B_n^{(g)}(u)$ holds. That is:
\begin{equation*}
    \cF_n^{(g)} := \bigcup_{u\in D_n} B_n^{(g)}(u).
\end{equation*}
Then $\cF_n^{(g)c}$ is the event that for each $u \in D_n$, at least one observation on the edges $(u,v)$ is not $X_n(v)^{-1}X_n(u)g$ or $(X_n(u)g)^{-1}X_n(v)$, that is, 
\begin{equation*}
    \cF_n^{(g)c} := \bigcap_{u\in D_n} B_n^{(g)}(u)^c.
\end{equation*}
Let $\{u_i\}_{i=1}^{|D_n|}$ enumerate distinct points from $D_n$; then for each $n \in \N$ it holds that 
\begin{equation*}
    \cF_n^{(g)c} \subset \bigcap_{i=1}^{|D_n|} B_n^{(g)}(u_i)^c.
\end{equation*}
Since $D_n$ is an independent set of vertices, $\{B_n^{(g)}(u_i)\}_{i=1}^{|D_n|}$ is collection of independent events. Then by monotonicity and this independence, we have
\begin{equation*}
    0 \leq \P(\cF_n^{(g)c}) \leq \prod_{i=1}^{|D_n|} \P(B_n^{(g)}(u_i)^c)  \le \left(1-\left(\frac{p_n}{|\fG|-1}\right)^{d}\right)^{|D_n|}.
\end{equation*}
Now we let $n$ tend to infinity. By hypothesis, we have $|D_n| \to \infty$, and by Lemma 4 (with $K = |\fG|-1 \ge 1 > 1/2$), we have 
\begin{equation*}
    \lim_{n \to \infty} \left(1-\left(\frac{p_n}{|\fG|-1}\right)^{d}\right)^{|D_n|} = 0.
\end{equation*}
Then by the Squeeze Theorem we have
\begin{equation*}
    \lim_{n \to \infty} \P(\cF_n^{(g)c}) = 0 \implies \lim_{n \to \infty} \P(\cF_n^{(g)}) =1.
\end{equation*}
Therefore, as $n$ goes to infinity, with probability tending to one, there exists a vertex $v \in V_n$ such that each observation $Y_n(u,v)$ or $Y_n(v,u)$ is $X_n(u)^{-1}X_n(v)g$ or $(X_n(v)g)^{-1}X_n(u)$ respectively. 
\end{proof}

\begin{PfofTheorem2}

Assume the hypotheses of Theorem 2. Let $g \in \fG \setminus \{e\}$. 
Let $B_n^{(g)}(u_0)$ be defined as in the proof of Lemma 5. Suppose for the moment that this event occurs.
Define $\bar{X}_n^{(g)}$ by the rule: 
\begin{equation*}
  \bar{X}_n^{(g)}(u) = \begin{cases} 
      X_n(u_0)g & \text{if } u = u_0\\
      X_n(u) & \text{otherwise}.
   \end{cases}
\end{equation*}
Its clear from the definition of $\pi_n'$ that:
\begin{equation*}
    \frac{\pi_n'(X_n\mid Y_n)}{\pi_n'(\bar{X}_n^{(g)}\mid Y_n)} = \frac{Q_n(Y_n\mid\psi(X_n))}{Q_n(Y_n\mid\psi(\bar{X}_n^{(g)}))}.
\end{equation*}
Denote the set of edges adjacent to $u_0$ by $\cN_n(u_0)$. Evidently, for all $(u,v) \notin \cN_n(u_0)$,
\begin{equation*}
    Q_n'(Y_n(u,v)\mid\psi(X_n)(u,v)) =  Q_n'(Y_n(u,v)\mid\psi(\bar{X}_n^{(g)})(u,v)), 
\end{equation*}
and therefore
\begin{align*}
\frac{Q_n(Y_n\mid\psi(X_n))}{Q_n(Y_n\mid\psi(\bar{X}_n^{(g)}))}   &  = \frac{\prod_{(u,v) \in \cN_n(u_0)} Q_n'(Y_n(u,v)\mid\psi(X_n)(u,v))}{\prod_{(u,v) \in \cN_n(u_0)} Q_n'(Y_n(u,v)\mid\psi(\bar{X}_n^{(g)})(u,v))} \\
  &  = \frac{\left(\frac{p_n}{|\fG|-1} \right)^{|\cN_n(u_0)|}}{(1-p_n)^{|\cN_n(u_0)|}} < 1. \\
    \end{align*}
From the inequality, it follows
\begin{equation*}
    \pi_n'(X_n|Y_n) < \pi_n'(\bar{X}_n^{(g)}\mid Y_n).
\end{equation*}
By Lemma 3, this implies that 
\begin{equation*}
    \psi(X_n) \ne \psi(\Phi_n(Y_n)).
\end{equation*}
Note that we have shown that if $B^{(g)}_n(u_0)$ occurs, then the MAP estimator is incorrect, i.e., $\psi(X_n) \ne \psi(\Phi_n(Y_n))$.

Now let $g \in \fG \setminus \{e\}$ be arbitrary. 
Let $\cF_n^{(g)}$ be defined as in the proof of Lemma $5$: 
\begin{equation*}
    \cF_n^{(g)} = \bigcup_{v \in V_n} B_n^{(g)}(v).
\end{equation*}
Note that $\{\psi(\Phi_n(Y_n)) = \psi(X_n)\} \subset \cF_n^{(g)c}$.
By the optimality property of the MAP estimator and the above inclusion, the following holds: for any estimator $T_n:\fG^{E_n} \to \fG^{V_n}$,
\begin{equation*}
   0 \le \P(\psi(T_n(Y_n)) = \psi(X_n)) \le \P(\psi(\Phi_n(Y_n)) = \psi(X_n)) \le \P(\cF_n^{(g)c}).
\end{equation*}
By Lemma 5 we have that 
\begin{equation*}
    \lim_{n \to \infty} \P(\cF_n^{(g)}) = 1,
    \end{equation*}
    and therefore
    \begin{equation*}
        \lim_{n \to \infty} \P(\cF_n^{(g)c}) = 0.
\end{equation*}
Therefore, by the Squeeze Theorem, the probability that any sequence of estimators $\{T_n\}$ on $\{\fQ_n\}$ is correct tends to zero as $n$ tends to infinity, and we conclude that $\{\fQ_n\}$ does not have the strong recovery property. 
\end{PfofTheorem2}
\section{Conclusion/Open Questions}
In conclusion, we have established both positive and negative results about strong recovery in the group synchronization problem. For our positive result, we showed that if the graph is the complete digraph on $n$ vertices and the noise is bounded away from $\sfrac{1}{|\fG|}$, then there exists a polynomial time algorithm which solves the synchronization problem with probability tending to one as $n$ goes to infinity. In our negative result, we have shown that if $G_n$ has an independent set $D_n \subset V_n$ whose size tends to infinity, and if the noise does not tend to zero too quickly, then error propagation ensures that strong recovery is impossible with probability tending to one as $n$ tends to infinity. We consider the following questions interesting avenues for future research. Can the result for sequences of complete digraphs be generalized to hold for a more broad class of graphs? In particular, perhaps the estimator of the proof of Theorem 1 can be applied to dense graphs or expander graphs?  

\bibliographystyle{abbrv}

\end{document}